\documentclass{tac}
\usepackage{xy}
\usepackage{MnSymbol}
\usepackage{amsfonts}

\usepackage{tikz}

\xyoption{all}

\def\mathcaldef#1{\expandafter\def\csname#1\endcsname{{\cal#1}}}

\def\"{``}

\def\qq{\quad\quad}
\def\qv{\qq ;\qq}

\def\iso{\,\cong\,}

\def\adj{\dashv}
\def\op{^{\rm op}}
\def\ov{\overline}
\def\ul{\underline}

\def\al{\alpha}

\def\dSet{\mathbb{S}\rm{et}}
\def\dSpan{\mathbb{S}\rm{pan}\,}
\def\dCat{\mathbb{C}\rm{at}}
\def\dSq{\mathbb{S}\rm{q}\,}
\def\dPb{\mathbb{P}\rm{b}\,}
\def\dA{\mathbb A}
\def\dB{\mathbb B}

\def\dQ{\mathbb Q\,}

\mathcaldef{A}
\mathcaldef{B}
\mathcaldef{C}
\mathcaldef{D}
\mathcaldef{S}

\mathbfdef{Set}
\mathbfdef{Set_f}
\mathbfdef{cMon}
\mathbfdef{Span}
\mathbfdef{Cat}
\mathbfdef{Prof}
\mathrmdef{id}

\def\pb{{\rm pb}}

\newtheorem{prop}{Proposition}

\let\pf\proof
\let\epf\endproof
\def\eq{\begin{equation}}
\def\eeq{\end{equation}}

\author{Claudio Pisani}
\address{via Saluzzo 67,\\ 10125 Torino, Italy.}
\title{Operads as double functors}
\copyrightyear{2022}
\eaddress{pisclau@yahoo.it}
\keywords{Double categories; operads; indexed categories; Beck-Chevalley condition}
\amsclass{18B10, 18C10, 18D30, 18D60, 18M60, 18M65, 18N10}

\begin{document}

\maketitle

\begin{abstract}

It is shown how double categories provide a direct abstract approach
to coloured operads; namely, product-preserving normal lax functors
$M:(\dPb \C)\op \to \dCat$ can be seen as generalized operads, 
the standard ones arising when $\C = \Set_f$.
In this context, generalized symmetric monoidal categories are considered,
in particular those arising from indexed categories with sums or products.
\end{abstract}


\section{Introduction}
\label{intro}

The role of double categories as a framework for abstract or formal
category theory has been promoted by several authors 
(notably by Bob Paré and Mike Shulman) and seems by now widely accepted.
Indeed, it is well known that many categorical notions and properties can be treated 
in any \"virtual equipment", playing the role of the double category $\dCat$
of functors and profunctors. 

Furthermore, \"all" generalized multicategories 
can be seen as instances of $T-$monoids in some virtual double category 
$\mathbb{D}$, where $T$ is a monad on $\mathbb{D}$
(see in particular \cite{leinster} and  \cite{shulman}).
While this approach has the advantage of bringing under the same umbrella 
many sorts of multicategories, at the same time it may result too general 
when one is interested in some specific cases.
Indeed, the relevant aspects are in a sense hidden in the monad $T$, 
and the monad itself often looks rather involved or unmanageable. 
This is the case in particular for one of the most commonly arising sort 
of multicategories, namely the symmetric ones, also known as (coloured) operads.

The thesis of the present paper is that the theory of double categories is actually
more directly involved in summarizing some aspects of (symmetric) multicategories 
than the above mentioned perspective could suggest.
In fact, classical coloured operads 
(in their non-skeletal form) are simply product-preserving (double lax) functors 
\eq
\label{int1}
M : (\dPb \Set_f)\op \to \dSet 
\eeq
to the double category of mappings and spans. 
Thus $(\dPb \Set_f)\op$, the horizontal dual of the double category of 
pullback squares in finite sets, is, in a sense, the \"double theory" for operads.
This can be seen as a generalization of the fact that product-preserving functors
\eq
\label{int2}
M : (\dPb \Set_f)\op \to \dSq\Set 
\eeq
are commutative monoids, where $\Set$ can be replaced by any 
finite-product category $\S$.

By the universal property of the module construction (see \cite{shulman}), 
which applied to $\dSet$ gives $\dCat$,
one can equivalently consider, instead of (\ref{int1}), normal functors
\eq 
\label{int3}
M : (\dPb \Set_f)\op \to \dCat 
\eeq
While the first perspective is more suitable to treat the fibered aspects, 
this one has technical advantages and is preferred in this work. 

Of course, the present approach points toward a notion of generalized operad 
rather different from the monad-based ones; 
namely we are led to consider normal product-preserving functors
\eq 
\label{int4}
M : (\dPb \C)\op \to \dCat 
\eeq
(or to $\dCat\,\S$) as the basic concept. 
We call such a functor a \"DF operad" (on $\C$), where the prefix \"DF" 
stands for \"double functor" and will accompany
also the various concepts generalized in this context. 

Thus, the horizontal part $M_h$ of a DF operad $M$ consists 
of an indexed category $M:\C\op\to\Cat$ 
(giving the $I$-\"families" $MI$ which are the domains of the \"operations")
such that 
\eq
\label{int5}
M(\Sigma I) \iso \Pi (MI)
\eeq
while the vertical part $M_v$ gives the operations themselves and their
compositions.

Among DF operads we isolate the \"DF (symmetric) monoidal categories" 
(in which the profunctors $M_vf$ are representable and the associated
Beck-Chevalley condition holds)
and the \"DF fibrations" (in which $M_vf$ is co-represented by $M_hf$).

Typically, DF monoidal categories arise from an indexed categories 
$M:\C\op\to\Cat$ satisfying (\ref{int5}) and having indexed sums or products;
that is, such that substitution (or reindexing) functors $f^* = M_h f$ 
have left or right adjoint 
\[
f_!\adj f^* \qv f^*\adj f_*
\]
satisfying the Beck-Chevalley law.
In the case of indexed sums, the corresponding DF monoidal categories
are also DF fibrations, since profunctors represented by $f_!$ are also
co-represented by $f^*$.

In the classical case $\C = \Set_f$, $M_h$ is the standard indexing 
$\A^I$ for a category $\A$. If $\A$ has finite sums or products,
the above mentioned DF monoidal categories are of course the
usual cocartesian or cartesian monoidal categories.
On the other and, if the category $\A$ has all small (co)products, we can
use $\Set$ in place of $\Set_f$ as the indexing category $\C$
and we still get a (genuinely generalized) DF operad.  
In this perspective, the role of the Beck-Chevalley condition is to assure that 
by taking isomorphism classes of a DF monoidal category one gets 
a (genuinely generalized) DF monoid 
\eq
M : (\dPb \C)\op \to \dSq\Set 
\eeq
(recall (\ref{int2}) and see section \ref{mon}).
This makes it precise the idea that also infinite categorical sums or products 
give a well defined algebraic structure.

The present work originates from \cite{pisani2}, whose emphasis was on the double
fibered form of operads (see also \cite{lambert}).


\section{Some remarks on double categories}
\label{dc}

We assume the basic notions about double categories and lax functors 
(see for instance \cite{pare}, \cite{pare2} and \cite{shulman2}).
In the present section we fix notations and highlight the facts that will be
useful in the sequel; some of them (especially those in section \ref{mackey})
 seem to be new, or at least not widely known.


\subsection{Notations and conventions}\qq
All double categories will be written with the first letter in a blackboard style.
Since a double category $\dA$ is a category in $\Cat$, it has an underlying
graph in $\Cat$ with source and target functors:
\eq 
\label{dc2}
s : A_1 \to A_0 \qv t : A_1 \to A_0 
\eeq
We adopt the convention that the arrows of $A_0$ are called \"horizontal arrows",
while the objects of $A_1$ are called \"vertical arrows" or \"proarrows"
and the arrows of $A_1$ are the \"cells" or \"squares". 
Thus a typical cell $\alpha:f\to g$ in $A_1$, with $s\alpha = k$ and $t\alpha = l$,
can be depicted as 
\eq
\label{dc1}
\xymatrix@R=1.3pc@C=1.3pc{
I \ar[dd]_f\ar[rr]^k  & & L \ar[dd]^g \\
& {\alpha} & \\
J  \ar[rr]_l                 & & K       }
\eeq

By a double category $\dA$ we intend a pseudo (or weak) one: 
the vertical or \"external" composition $\circ : A_1\times_{A_0}A_1\to A_1$
is defined up to isomorphisms, as in bicategories. 
By a functor $F:\dA \to \dB$ we intend a lax one (see \cite{pare} or \cite{shulman2}).
A functor is normal if identities are strictly preserved.
We denote by $F_h$ and $F_v$ the horizontal and the vertical components of $F$:
\[
\xymatrix@R=1.3pc@C=1.3pc{
I \ar[dd]_f\ar[rr]^k  & & L \ar[dd]^g \\
& {\alpha} & \\
J  \ar[rr]_l                 & & K       }
\xymatrix@R=2pc@C=2pc{
 \\ & \ar@{|->}[r]^F &     } \qq
\xymatrix@R=1.3pc@C=1.3pc{
FI \ar[dd]_{F_vf}\ar[rr]^{F_h k}  & & FL \ar[dd]^{F_vg} \\
& {F\alpha} & \\
FJ  \ar[rr]_{F_hl}                & & FK       }
\]
By the opposite double category $\dA\op$ we mean the horizontal opposite,
obtained by taking the opposite of $A_0$ and of $A_1$ 
(but not reversing $s$ and $t$):
\[ s\op : A_1\op \to A_0\op \qv t\op : A_1\op \to A_0\op \]
Thus, a cell (\ref{dc1}) in $\dA$ corresponds to the following cell in $\dA\op$:
\[
\xymatrix@R=1.3pc@C=1.3pc{
I \ar[dd]_f  & & L \ar[dd]^g\ar[ll]_k \\
& {\alpha} & \\
J                  & & \ar[ll]^l K       }
\]


Following Paré, $\dSet = \dSpan\Set$ is the double category of mappings 
and spans and $\dCat$ is the double category of functors and profunctors. 
More generally, for a category $\S$ with pullbacks, $\dCat\,\S$ 
is the double category of internal categories, functors and profunctors, 
obtained by applying the module construction to the double category $\dSpan\S$.
For a category $\C$, $\dSq\C$ is the double category of commutative squares in $\C$, 
and $\dPb\C$ is the double category of pullback squares in $\C$.


\subsection{Double categories with products}
\label{pro}
\qq
Following \cite{pare2} we say that a double category $\dA$ has
(finite) products if both $A_0$ and $A_1$ have (finite) products
and the source and target functors (\ref{dc2}) preserve them.

So, for any pair of proarrows $f:I\to J$ and $g:K\to L$ we have a diagram
\eq   
\label{dc4}
\xymatrix@R=1.3pc@C=1.3pc{
I \ar[dd]_f  && I\times K \ar[dd]\ar[ll]_{p_1}\ar[rr]^{p_2} && K \ar[dd]^g \\
& \pi_1 && \pi_2 \\
J     && J\times L \ar[ll]^{q_1} \ar[rr]_{q_2}     && L       }
\eeq
encompassing a product diagram in $A_1$ along with its source and target 
diagrams in $A_0$ which also are product diagrams.
\begin{remark}
There are also stronger notions of product in a double category;
this one is best suited to our needs.
\end{remark}
Similarly, if $\dA$ has sums we have diagrams
\[
\xymatrix@R=1.3pc@C=1.3pc{
I \ar[rr]^{i_1}\ar[dd]_f  && I+K \ar[dd] && K \ar[dd]^g\ar[ll]_{i_2}  \\
& \iota_1 && \iota_2 \\
J \ar[rr]_{j_1}       && J+L   && L \ar[ll]^{j_2}       }
\]
Clearly, $\dA$ has sums if and only if $\dA\op$ has products.
\begin{prop}
If $\C$ has (finite) products, then $\dA = \dSq\C$ also as (finite) products;
the same holds true for sums.
\end{prop}
\pf
Indeed, in $\dA = \dSq\C$ the central vertical arrow of diagram (\ref{dc4}) 
is necessarily
$f\times g$, the unique arrow making the diagram commute in $\C$,
and it gives the product in $A_1$ (the arrow category of $\C$).
The same proof holds for sums. Note that we cannot refer naively 
to duality, since $\dSq(\C\op) \neq (\dSq\C)\op$; rather,
$\dSq(\C\op)$ is the vertical dual of $(\dSq\C)\op$ ($s$ and $t$ are also reversed).
But taking the vertical dual does not affect the existence of products or sums.
\epf
\begin{remark}
Actually, the proposition is true for all those double categories $\dA$ 
constructed from a bipointed category $\sigma,\tau:1\to D$ with
$A_0 = \C$, $A_1 = \C^D$ 
and with $s,t : A_1\to A_0$ induced by $\sigma$ and $\tau$.
If $\C$ has (finite) products, then $\dA$ has (finite) products, since
products in $\C^D$ are computed pointwise.
If $D$ is the \"free arrow" 
$\xymatrix@R=1.3pc@C=1.3pc{\bullet \ar[r]  & \bullet}$ we get $\dSq\C$,
while if $D$ is the \"free span"
$\xymatrix@R=1.3pc@C=1.3pc{\bullet  & \bullet \ar[r]\ar[l] & \bullet }$,
we get $\dSpan\C$.
Note that while the vertical structure for $\dSq\C$ is induced by the pushout 
(in $\Cat$) of the free arrow with itself, the same is not true for $\dSpan\C$.
It would be interesting to find general procedures to construct the vertical structure
for $A_1 = \C^D$ for some sorts of bipointed categories $D$.
\end{remark}

We say that a functor $F:\dA \to \dB$ preserves products 
(or that it is product-preserving) if both its components
$F_h : A_0 \to B_0$ and $F_v : A_1 \to B_1$ preserve products.
For instance, both the \"inclusion" functors 
\[ \dSq\C \to \dSpan\C \] 
preserve products, as can be easily seen directly; 
alternatively, one may note that these functor are induced by each of the two 
bipoint-preserving functors from the free span to the free arrow. 

\begin{prop}
\label{prdc1}
If the horizontal component $F_h:A_0\to \C$ of a functor
$F:\dA \to \dSq\C$ preserves products, then $F$ itself preserves products.
\end{prop}
\pf
The effect of $F$ on a product diagram is 
\[
\xymatrix@R=1pc@C=1pc{
I \ar[dd]_f  && I\times K \ar[dd]\ar[ll]_{p_1}\ar[rr]^{p_2} && K \ar[dd]^g \\
& \pi_1 && \pi_2 \\
J     && J\times L \ar[ll]^{q_1} \ar[rr]_{q_2}     && L       }
\xymatrix@R=2pc@C=2pc{
 \\ & \ar@{|->}[r]^F &     } \qq
\xymatrix@R=1pc@C=1pc{
FI \ar[dd]_{F_v f}  && F(I\times K) \ar[dd]\ar[ll]_{F_h p_1}\ar[rr]^{F_h p_2} && FK \ar[dd]^{F_vg} \\
& F\pi_1 && F\pi_2 \\
FJ     && F(J\times L) \ar[ll]^{F_h q_1} \ar[rr]_{F_h q_2}     && FL       }
\]
By the hypothesis, $F(I\times K) = FI\times FK$ and $F(J\times L) = FJ\times FL$,
that is, the horizontal diagrams on the right are products in $\C$.
Since diagrams in $\dSq\C$ are (commutative) diagrams in $\C$, the vertical
central arrow is forced to be $F_v f \times F_v g$.
\epf


\subsection{The double category $\dPb\C$ and Mackey functors}\qq
\label{mackey}
Recall that an extensive category is a category $\C$ with finite sums that 
interact nicely with pullbacks (see \cite{carboni}):
given a commutative diagram in $\C$ 
\[
\xymatrix@R=1.3pc@C=1.3pc{
I \ar[rr]\ar[dd]  && R \ar[dd] && K \ar[dd]\ar[ll]  \\
& &&  \\
J \ar[rr]_{i_1}       && J+L   && L \ar[ll]^{i_2}       }
\]
where the bottom row is a sum, 
the top row is also a sum if and only if the two squares are pullbacks.


\begin{prop}
If $\C$ is extensive, then $\dPb\C$ has sums and the inclusion $\dPb\C \to \dSq\C$ preserves them.
\end{prop}
\pf
By the extensivity of $\C$, the following is a sum diagram both in $\dSq\C$
and in $\dPb\C$:
\[
\xymatrix@R=1.3pc@C=1.3pc{
I \ar[rr]^{i_1}\ar[dd]_f  && I+K \ar[dd]|{f+g} && K \ar[dd]^g\ar[ll]_{i_2}  \\
& \pb && \pb \\
J \ar[rr]_{j_1}       && J+L   && L \ar[ll]^{j_2}       }
\]
\epf
Given two categories $\C$ and $\D$, a Mackey functors $\C\to\D$ 
in the sense of \cite{lidner} consists of two functors which coincide on objects
\[ F_!:\C\to\D \qv F^*:\C\op\to\D \]
with $F^*$ preserving finite products 
and such that for any pullback in $\C$
the corresponding right hand square below commutes:
\[
\xymatrix@R=1.3pc@C=1.3pc{
I \ar[dd]_f\ar[rr]^k  & & L \ar[dd]^g \\
& \pb & \\
J  \ar[rr]_l                 & & K       }
\qq\qq
\xymatrix@R=1.8pc@C=1.8pc{
FI \ar[dd]_{F_!f}  & & FL \ar[ll]_{F^*k}\ar[dd]^{F_!g} \\
&  & \\
FJ                  & & \ar[ll]^{F^*l} FK       }
\]

The following proposition is then an easy consequence of the definitions.
\begin{prop}
\label{prodc2}
For an extensive category $\C$, the Mackey functors $\C\to\D$ 
correspond to the finite-product-preserving functors $(\dPb\C)\op \to \dSq\D$.
\end{prop}
\pf
The correspondence is obviously given by $F_! \longleftrightarrow F_v$ 
and $F^* \longleftrightarrow F_h$. 
By proposition \ref{prdc1}, preservation of products $\C\op \to \D$ is sufficient
to give preservation of products $(\dPb\C)\op \to \dSq\D$.
\epf


\subsection{The double category $\dCat$}\qq
Recall that $\dCat$ is the double category of functors and profunctors.
Following \cite{joyal} and  \cite{pare}, a profunctor $\Phi:\C \to \D$ is a functor $\C\op\times\D \to \Set$,
rather than $\D\op\times\C \to \Set$ as many authors prefer.
One reason for this choice is that profunctors $\C \to \D$
can be seen as categories over the free arrow category 
$\xymatrix@R=1.3pc@C=1.3pc{0 \ar[r]  & 1}$ (the \"barrels" of \cite{joyal}),
and in this way the domain and the codomain of the profunctor
are projected respectively on the domain and the codomain of the free arrow.

A cell in $\dCat$ 
\[
\xymatrix@R=1.3pc@C=1.3pc{
\A \ar[dd]_\Phi\ar[rr]^F  & & \C \ar[dd]^\Psi \\
& {\alpha} & \\
\B  \ar[rr]_G                 & & \D       }
\]
is a natural transformation $\al : \Phi \to \Psi(F - , G -)$, that is a family of mappings
\[ \al_{X,Y} : \Phi(X,Y) \to \Psi(FX , GY) \qv X\in\C , \,\, Y\in\D \]
satisfying the usual conditions.
Equivalently, it is a morphism of barrels, that is of categories over the free arrow.

A profunctor $\Phi:\C\to\D$ is \"representable" if it isomorphic to 
$\Phi_F = \D(F-,-)$ for a functor $F:\C\to\D$.
As a barrel, $\Phi$ is representable if it is an opfibration over the free arrow.
The assignment $F \mapsto \Phi_F$ defines a functor $\dSq\Cat \to \dCat$.
For any square in $\Cat$, we get a cell $\al$ in $\dCat$
\[
\xymatrix@R=1.3pc@C=1.3pc{
\A \ar[dd]_F\ar[rr]^H  & & \C \ar[dd]^G \\
&  & \\
\B  \ar[rr]_L                 & & \D       }
\xymatrix@R=2pc@C=2pc{
 \\ & \ar@{|->}[r] &     } \qq
\xymatrix@R=1.3pc@C=1.3pc{
\A \ar[dd]_{\Phi_F}\ar[rr]^H  & & \C \ar[dd]^{\Phi_G} \\
& \al & \\
\B  \ar[rr]_L                 & & \D       }
\]
where the $\al_{X,Y} : \Phi_F(X,Y) \to \Phi_G(HX , LY)$, that is 
\[ \al_{X,Y} : \B(FX,Y) \to \D(GHX , LY) = \D(LFX , LY) \]
are given by the action of $L$ on the arrows of $\B$.

\begin{prop}
$\dCat$ has products, and the \"inclusion"
$\dSq\Cat \to \dCat$ preserves products.
\end{prop}
\pf
A product diagram in $\dCat$ is 
\[
\xymatrix@R=1.3pc@C=1.3pc{
\A \ar[dd]_\Phi  && \A\times\C \ar[dd]|{\Phi\times\Psi}\ar[ll]_{p_1}\ar[rr]^{p_2} 
&& \C \ar[dd]^\Psi \\
& \pi_1 && \pi_2 \\
\B     && \B\times\D \ar[ll]^{q_1} \ar[rr]_{q_2}     && \D       }
\]
where $\Phi\times\Psi: \A\times\C \to \B\times\D$ is defined by 
\[ (\Phi\times\Psi) ((X,X'),(Y,Y')) = \Phi (X,Y) \times \Psi(X',Y') \]
That this is indeed a product in $\dCat$ is easily checked directly.
Alternatively note that, for $\dA = \dCat$, $A_1$ is the category of barrels and
the product of $\Phi$ and $\Psi$ over the free arrow corresponds to the above one.
Similarly, the fact that $\dSq\Cat \to \dCat$ preserves products is easily 
checked directly. Alternatively note that the product of two representable
profuncors (as barrels which are opfibrations) is again representable,
by the product functor.
\epf


\section{DF operads}
\label{go}

In this section we study our version of generalized operad, namely
\"double functor operads" or briefly DF operads.
By an operad we intend a coloured operad, that is a symmetric multicategory.
In order to introduce the idea, and also because some results 
will be used in the sequel, we begin by presenting the relative notion
of generalized commutative monoid.


\subsection{DF monoids}
\label{mon} \qq
Given an extensive category $\C$ with pullbacks and a category $\S$
with products, a {\em DF monoid} (on $\C$ in $\S$) is a 
product-preserving functor
\[ F : (\dPb\C)\op \to \dSq\S \]
By Proposition \ref{prodc2}, DF monoids correspond to Mackey functors
$\C \to \S$ in the sense of \cite{lidner}, except for the fact that we
require preservation of all products, not only of the finite ones.
Thus, they consist of functors
\[ F_!:\C\to\S \qv F^*:\C\op\to\S \]
coinciding on objects, with $F^*$ preserving products 
and such that for any pullback in $\C$
the corresponding right hand square below commutes:
\eq
\label{go1}
\xymatrix@R=1.3pc@C=1.3pc{
I \ar[dd]_f\ar[rr]^k  & & L \ar[dd]^g \\
& \pb & \\
J  \ar[rr]_l                 & & K       }
\qq\qq
\xymatrix@R=1.8pc@C=1.8pc{
FI \ar[dd]_{F_!f}  & & FL \ar[ll]_{F^*k}\ar[dd]^{F_!g} \\
&  & \\
FJ                  & & \ar[ll]^{F^*l} FK       }
\eeq
In order to justify our terminology, we have:
\begin{prop}
DF monoids on $\Set_f$ in $\S$ coincide with commutative monoids in $\S$.
\end{prop}
\pf
We give two proofs. The first one appeals to the well-known result in \cite{lidner}:
if $\C$ is extensive with pullbacks, Mackey functors $\C \to \S$ correspond to 
finite-product-preserving functors $\Span(\C) \to \S$.
But $\Span(\Set_f)$ is essentially the category of matrices valued in natural numbers,
which is well known to be the Lawvere theory for commutative monoids.

The second proof is more direct and maybe more instructive.
We assume for simplicity $\S = \Set$, but the proof can be easily 
adapted to the more general case.
First, since $F^*:\Set_f\op \to \Set$ preserves products and $\Set_f\op$ 
is the free finitely complete category on the terminal category, $F^*$
is isomorphic to the family functor on an object $M\in\Set$;
thus we can assume $F^* I = M^I$, and we rewrite the right hand square 
in (\ref{go1}) as
\[
\xymatrix@R=1.8pc@C=1.8pc{
M^I \ar[dd]_{f_!}  & & M^L \ar[ll]_{k^*}\ar[dd]^{g_!} \\
&  & \\
M^J                  & & \ar[ll]^{l^*} M^K       }
\]
where we put $f^* = F^*f$ and $f_! = F_!f$.
Now, given a family $x_i,\, i\in I$ in $M^I$, we define its \"product"
$m_I(x_i)$ as the element $f_!(x_i)$ of $M = M^1$, where $f:I\to 1$ in $\Set_f$.
The product of a family is stable with respect to bijective reindexing $k^*$, 
since the pullback on the left induces a square on the right
\[
\xymatrix@R=1.3pc@C=1.3pc{
I \ar[dd]_f\ar[rr]^k  & & L \ar[dd]^g \\
& \pb & \\
1  \ar[rr]_\id                & & 1       }
\qq\qq
\xymatrix@R=1.8pc@C=1.8pc{
M^I \ar[dd]_{m_I}  & & M^L \ar[ll]_{k^*}\ar[dd]^{m_L} \\
&  & \\
M                  & & \ar[ll]^\id M       }
\]
Thus, for $a,b\in M$, the element $m(a,b)\in M$ is not ambiguous 
and $m(a,b) = m(b,a)$; the same holds for, say, $m(a,b,c)$.
To show that the operation is associative, consider any three element set,
say $I = \{i,j,k\}$ and any two element set, say $J = \{r,s\}$.
Then, since $f:I\to 1$ factors through $J$ as $hg$ where $g:i,j\mapsto r$
and $g:k\mapsto s$, we have the diagram on the left, where the rows are sums, 
which corresponds to the diagram on the right:
\[
\xymatrix@R=1.3pc@C=1.3pc{
S \ar[dd]\ar[rr]  & & I \ar[dd]^g & & 1 \ar[dd]^\id\ar[ll] \\
& \pb & & \pb \\
1 \ar[rr]^r  & & J \ar[dd]^h & & 1 \ar[ll]_s \\
&  & \\
               & & 1       }
\qq\qq
\xymatrix@R=1.4pc@C=1.4pc{
M^S \ar[dd]_m  & & M^I \ar[ll]\ar[dd]^{g_!} & & M \ar[dd]\ar[ll] \\
&  & \\
M  & & M^J \ar[ll] \ar[dd]^m & & M \ar[ll] \\
&  & \\
             & & M       }
\]
Now, since a DF monoid $(\dPb\C)\op \to \dSq\S$ preserves products
by definition, $g_! = m\times\id$.
(In general, the product preservation axiom assures that 
the value of $(-)_!$ on an arbitrary mapping depends only from its value 
on the mappings with a terminal codomain.)
Since by functoriality $m g_! = m_I$, we have
\[
m(m(a,b),c) = m((m\times\id)(a,b,c)) = (m g_!)(a,b,c) = m_I(a,b,c) 
\]
Symmetrically, one has $m(a,m(b,c)) = m_I(a,b,c)$ as well, 
which proves associativity.
Along the same lines, one can verify the identity laws.

In the other direction, given a commutative monoid $M$ and a mapping 
$f:I\to J$ in $\Set_f$, it is obvious how to exploit the multiplication of $M$
in order to obtain a mapping $f_!:M^I\to M^J$ and it is easy to verify that
one so gets a DF monoid.
\epf


\subsection{DF operads}
\label{op} \qq
Given an extensive category $\C$ with pullbacks and a category $\S$
with pullbacks, a {\em DF operad} (on $\C$ in $\S$) is a 
product-preserving normal functor
\[ M : (\dPb\C)\op \to \dCat\S \]
In the sequel, we just consider the case $\S = \Set$, 
but the general case can be easily treated. 
Thus, a DF operad consists of a functor
$M_h:\C\op\to\Cat$ preserving products and a normal lax functor $M_v:\C\to\Prof$
coinciding with $M_h$ on objects and such that for any pullback in $\C$ 
there is a corresponding cell in $\dCat$:
\eq
\label{go3}
\xymatrix@R=1.3pc@C=1.3pc{
I \ar[dd]_f\ar[rr]^k  & & L \ar[dd]^g \\
& \pb & \\
J  \ar[rr]_l                 & & K       }
\qq\qq
\xymatrix@R=1.6pc@C=1.6pc{
MI \ar[dd]_{M_vf}  & & ML \ar[ll]_{M_hk}\ar[dd]^{M_vg} \\
& \al & \\
MJ                  & & \ar[ll]^{M_hl} MK       }
\eeq
The assignment should satisfy the axioms for lax double functors and should also
preserve products in the sense of section \ref{pro}.
In order to justify our terminology, we are going to show that DF operads
on $\Set_f$ are essentially classical coloured operads 
(that is symmetric multicategories) in their non-skeletal form.
The need for a non-skeletal form of operads was already stressed in
\cite{leinster}, where a version of them (named \"fat symmetric multicategories")
is given. We hope to convince the reader that the present version is  
very natural.
\begin{prop}
\label{prop}
DF operads on $\Set_f$ give the non-skeletal version of operads.
\end{prop}
\pf (Sketch.)
Let $\bf N$ be the usual skeleton of $\Set_f$, with objects ${\bf n} = \{1, 2, \dots, n\}$,
and suppose that a bijection is $b_I:I\to{\bf n}$ is given, for every $I\in\Set_f$.
In one direction, given a classical operad $M$ we get a DF operad $\ov M$ 
in the following way: $\ov MI$ is $M_0^I$ where $M_0$ is the underlying category
of $M$, and for a mapping $l:J\to K$ in $\Set_f$, $\ov M_hl:\ov MK\to\ov MJ$
is the standard reindexing functor. 
If $f:I\to J$ and $J$ is terminal in $\Set_f$, we get the profunctor 
$\Phi_f = \ov M_v f:M_0^I\to M$ by posing $\Phi_f(X;A) = M(A_1,\dots,A_n ; A)$
(where we use $b_I$ to reindex the $I$-family $X = A_i, i\in I$).
For a general $f:I\to J$ in $\Set_f$, we get $\Phi_f = \ov M_v f:M_0^I\to M_0^J$ 
by posing $\Phi_f(X;Y) = \Pi_j\Phi_{f_j}(X_j;B_j)$, where $Y = B_j, j\in J$,
and $f_j$ and $X_j$ are obtained by restricting $f$ and $X$ to $f^{-1}j\subseteq I$.
Laxity cells $\Phi_g \Phi_f \to \Phi_{gf}$ are of course given by composition
in $M$ and, for a pullback square in $\Set_f$, the corresponding cell $\al$ 
in (\ref{go3}) is obtained with the aid of the action of bijective mappings $\bf n \to n$
on the arrows of $M$, mediated by the $b_I$.

In the other direction, given a DF operad $M$ on $\Set_f$, we define 
a classical operad $\ul M$ which has $M = M1$ as objects set and
with $\ul M(A_1,\dots,A_n ; A) = (M_vf)(X,A)$, where $X = A_i, i\in\bf n$
and $f:\bf n\to 1$.
Composition in $\ul M$ is given by the laxity cells and actions of bijective 
mappings $k:\bf n \to n$ is given by the cell $\al$ below.
\[
\xymatrix@R=1.5pc@C=1.5pc{
{\bf n} \ar[dd]\ar[rr]^k  & & {\bf n} \ar[dd] \\
& \pb & \\
1  \ar[rr]_\id                & & 1       }
\qq\qq\qq
\xymatrix@R=1.6pc@C=1.6pc{
M{\bf n} \ar[dd]  & & M{\bf n} \ar[ll]_{M_hk}\ar[dd] \\
& \al & \\
M                  & & \ar[ll]^\id M       }
\]
The fact that DF operads $M : (\dPb\Set_f)\op \to \dCat$ preserve products,
assures that value of $M_v$ on arbitrary mappings $f:I\to J$ 
is determined by their value on mappings with a terminal codomain.
Indeed, $f = \Sigma_j f_j$ in $\dPb\Set_f$, so that 
$f = \Pi_j f_j$ in $(\dPb\Set_f)\op$ and $M_vf = \Pi_j M_vf_j$ in $\dCat$.
It follows that the DF operads $M$ and $\ov{(\ul M)}$ are essentially the same.
This complete our sketch of the proof.
\epf
\begin{remark}
Note that among the laws implicit in the lax functoriality of $M$,
there is the one regarding composition of a family of arrows with an arrow
acted over by a bijective mapping $k:\bf n \to n$:
composing along the vertical left side of the right hand 
diagram below gives the same result as composing along the vertical 
right side and then acting on the composition by the bijective mapping $t$.
Such a $t$, obtained by composing the pullbacks on the left, is exactly
the mapping which requires a rather involved explicit description in 
the usual definitions of symmetric multicategories.
\[
\xymatrix@R=1.6pc@C=1.6pc{
{\bf m} \ar[dd]\ar[rr]^t  & & {\bf m} \ar[dd] \\
& \pb & \\
{\bf n} \ar[dd]\ar[rr]^k  & & {\bf n} \ar[dd] \\
& \pb & \\
1  \ar[rr]_\id                & & 1       }
\qq\qq\qq
\xymatrix@R=1.6pc@C=1.6pc{
M{\bf m} \ar[dd]  & & M{\bf m} \ar[ll]_{M_ht}\ar[dd] \\
& \beta & \\
M{\bf n} \ar[dd]  & & M{\bf n} \ar[ll]_{M_hk}\ar[dd] \\
& \al & \\
M                  & & \ar[ll]^\id M       }
\]
\end{remark}

\begin{remark}
Commutative monoids (in $\Set$) can be seen as discrete
symmetric multicategories $M$, that are those such that the functor
$M\to 1$ is a discrete opfibration of multicategories. 
Explicitly, this means that for any $A_1,\dots,A_n \in M$ there is exactly
one arrow $A_1,\dots,A_n \to A$ out of them, giving their product $A\in M$.
Translating this particular case in the present framework, 
we note that the categories $MI$ are then discrete and the profunctors 
$\Phi_f$ are given, for $f:I\to J$, $X = (A_i)_{i\in I}$ and $Y = (B_j)_{j\in J}$, 
by $\Phi_f(X,Y) = 1$ if $B_j$ is the product (in the above sense) of 
$(A_i)_{i\in f^{-1}j}$, for any $j\in J$, and $\Phi_f(X,Y) = 0$ otherwise.
In other terms, the $M_vf = \Phi_f$ are in fact mappings assigning to
any $X\in MI = M^I$ the unique $Y$ for which $\Phi_f(X,Y) = 1$.
Furthermore, for this particular sort of categories (namely the discrete ones) 
and profunctors (which are represented by mappings), the cells in $\dCat$
are cells in $\dSq\Set$.

Thus, as a particular case of proposition \ref{prop}, 
we find again the characterization of commutative
monoids given in section \ref{mon}.
\end{remark}

\begin{remark}
As one may expect, the natural notion of morphism of lax functors
$F,G:\dA \to \dB$, which is a double natural transformation $F \to G$
(see \cite{pare} or \cite{shulman2}),
gives the right notion of morphism of DF operads $M,N : (\dPb\C)\op \to \dCat$,
so that we have the category ${\bf Op}_\C$.
A sharper form of proposition \ref{prop} should then yield an equivalence 
between the category of classical operads and ${\bf Op}_{\Set_f}$,
as in the analogous result in \cite{leinster} on fat symmetric multicategories.
\end{remark}

\begin{remark}
In the case $\C = \Set_f$, exponentiable DF operads $M : (\dPb\C)\op \to \dCat$
are those such that $M$ is a pseudo functor, that is laxity cells are isos.
Indeed, these are the promonoidal symmetric multicategories, which
coincide with the exponentiable ones (see \cite{pisani}).
It seems likely that this characterization of exponentiabilty holds for any base $\C$.
\end{remark}


\subsection{DF monoidal categories}
\label{mc} \qq
A DF operad $M : (\dPb\C)\op \to \dCat$ is a {\em DF monoidal category} 
if the following conditions are satisfied:
\begin{enumerate}
\item
the profunctors $\Phi_f = M_vf:MI \to MJ$ are representable:

for any $f:I\to J$ in $\C$, there is a functor $f_!:MI\to MJ$
such that
\[ \Phi_f(X,Y) \iso MJ(f_!X,Y) \]
\item
the Beck-Chevalley condition holds: 

for any pullback square in $\C$, as the left hand one below,
the square on the right commutes up to isomorphisms:
\eq
\label{B-C}
\xymatrix@R=1.3pc@C=1.3pc{
I \ar[dd]_f\ar[rr]^k  & & L \ar[dd]^g \\
& \pb & \\
J  \ar[rr]_l                 & & K       }
\qq\qq\qq
\xymatrix@R=1.7pc@C=1.6pc{
MI \ar[dd]_{f_!}  & & ML \ar[ll]_{k^*}\ar[dd]^{g_!} \\
&  & \\
MJ                  & & MK \ar[ll]^{l^*}        }
\eeq
\end{enumerate}
Of course, for $\C = \Set_f$, we get the usual symmetric monoidal categories,
given in the universal form of representable symmetric multicategories
(see \cite{hermida}, \cite{leinster} and \cite{pisani}).

In particular, $f_!$ can be given by sums or products for 
the horizontal indexed category
\[ 
I\mapsto MI \qv  f\mapsto f^*=M_hf
\]
as mentioned in the introduction, yielding cartesian or cocoartesian
DF monoidal categories.
\begin{prop}
By taking isomorphism classes of a DF monoidal category $M$, 
one gets a DF monoid $|M| : (\dPb\C)\op \to \dSq\Set$.
Namely, $|M|I = |MI|$, $|M|_hf = |f^*|$ and $|M|_vf = |f_!|$
\end{prop}
\pf
The Beck-Chevalley condition (\ref{B-C}) ensures that the right hand square
\[
\xymatrix@R=1.3pc@C=1.3pc{
I \ar[dd]_f\ar[rr]^k  & & L \ar[dd]^g \\
& \pb & \\
J  \ar[rr]_l                 & & K       }
\qq\qq\qq
\xymatrix@R=1.6pc@C=1.6pc{
|MI| \ar[dd]_{|f_!|}  & & |ML| \ar[ll]_{|k^*|}\ar[dd]^{|g_!|} \\
&  & \\
|MJ|                  & & |MK| \ar[ll]^{|l^*|}        }
\]
commutes in $\Set$. 
\epf

\begin{remark}
By proposition \ref{prop}, DF monoids in $\Cat$
\[ M : (\dPb\Set_f)\op \to \dSq\Cat \]
correspond to commutative monoids in $\Cat$, that is strict 
symmetric monoidal categories.
It would be nice if general symmetric monoidal categories
could be captured by functors
\[ M' : (\dPb\Set_f)\op \to \dSq^*\Cat \]
where $\dSq^*\Cat$ is some kind of weak or pseudo version of $\dSq\Cat$.
For instance, one could check wether something like the subcategory $\dQ^*\Cat$ 
of the of the double category $\dQ\Cat$ of quintets in $\Cat$, formed by those 
cells which are natural isomorphisms, serves the purpose. 
In this case, one could simply say that DF monoidal categories are those
DF operads $M : (\dPb\C)\op \to \dCat$ which factors through a functor
$M' : (\dPb\C)\op \to \dSq^*\Cat$.
\end{remark}


\subsection{DF fibrations}
\label{fib} \qq
A {\em DF fibration} is a DF operad $M : (\dPb\C)\op \to \dCat$
such that the profunctors $\Phi_f = M_vf:MI \to MJ$ are co-represented
by $f^* = M_hf$, for any $f:I\to J$ in $\C$:
\[ \Phi_f(X,Y) \iso MI(X,f^*Y) \]
Thus, as expected, a DF fibration is determined by its horizontal component,
that is the indexed category 
\[ 
I\mapsto MI \qv  f\mapsto f^*
\]
The vertical component, which as a normal lax functor $\C\to\Prof$ corresponds 
to a category over $\C$, embodies the associated fibration.


%
%


\begin{refs}

\bibitem[Carboni et al.]{carboni} A. Carboni, S. Lack, R. F. C. Walters (1993), 
Introduction to extensive and distributive categories, {\em J. Pure Appl. Algebra} 
{\bf 84}, 145-158.

\bibitem[Cruttwell \& Shulman, 2010]{shulman} G. Cruttwell and M. Shulman (2010), 
A unified framework for generalized multicategories, {\em Theory and Appl. of Cat.} 
{\bf 24}, 580-655.

\bibitem[Hermida, 2000]{hermida} C. Hermida (2000), Representable multicategories, 
{\em Advances in Math.}, {\bf 151}, 164-225.

\bibitem[Joyal, 2022]{joyal} A. Joyal (last revised 2022), Distributors and barrels,
available on Joyal's CatLab in nlab.

\bibitem[Lambert, 2021]{lambert} M. Lambert (2021), Discrete double fibrations, 
preprint available on arXiv.org.

\bibitem[Leinster, 2003]{leinster} T. Leinster (2003), {\em Higher operads, higher categories}, Cambridge University Press, math.CT/0305049.

\bibitem[Lidner, 1976]{lidner}
H. Lindner (1976), A remark on Mackey-functors, {\em Manuscripta Math.} 
{\bf 18}, 273–278. 

\bibitem[Paré, 2011]{pare} R. Paré (2011), Yoneda theory for double categories, 
{\em Theory and Appl. of Cat.} {\bf 17}, 436-489.

\bibitem[Paré, 2018]{pare2} R. Paré (2018), Double Categories - 
The best thing since slice categories, slides for the 2018 FMCS workshop
at Dalhousie University, available on Bob Paré home page.  

\bibitem[Pisani, 2014]{pisani} C. Pisani (2014), Sequential multicategories, 
{\em Theory and Appl. of Cat.} {\bf 29}, 496-541.

\bibitem[Pisani, 2022]{pisani2} C. Pisani (2022), Fibered multicategory theory, 
preprint available on arXiv.org.

\bibitem[Shulman, 2008]{shulman2} M. Shulman (2008), 
Framed bicategories and monoidal fibrations, {\em Theory and Appl. of Cat.} 
{\bf 20}, 650-738.




\end{refs}

\end{document}